\documentclass[11pt, reqno]{amsart}
\usepackage[utf8]{inputenc}
\usepackage{amsmath}
\usepackage{amsthm}
\usepackage{amsfonts}
\usepackage{amssymb}
\usepackage{blindtext,amsmath,amsthm,amsfonts,amssymb, textcomp, mathrsfs, mathtools}
\DeclarePairedDelimiter \floor {\lfloor} {\rfloor}

\usepackage[top=35mm,bottom=35mm,left=35mm,right=35mm]{geometry}
\newtheorem{theorem}{Theorem}[section]
\newtheorem{conjecture}{Conjecture}
\newtheorem{corollary}{Corollary}

\newtheorem*{theorem*}{Theorem}
\newtheorem*{remark*}{Remark}
\newtheorem*{problem*}{Problem}
\newtheorem*{conjecture*}{Conjecture}
\newtheorem{lemma}[theorem]{Lemma}

\begin{document}
\title[On a conjecture of Livingston]{On a conjecture of Livingston}

\author[Siddhi Pathak]{Siddhi Pathak}

\address{Department of Mathematics and Statistics, Queen's University, Kingston, Canada, ON K7L 3N6.}
\email{siddhi@mast.queensu.ca}



\begin{abstract}
In an attempt to resolve a folklore conjecture of Erd\"os regarding the non-vanishing at $s=1$ of the $L$-series
attached to a periodic arithmetical function with period $q$ and values in $\{ -1, 1\} $, Livingston conjectured the $\bar{\mathbb{Q}}$ - linear independence of logarithms of certain algebraic numbers. In this paper, we disprove Livingston's conjecture for composite $q \geq 4$, highlighting that a new approach is required to settle Erd\"os’s conjecture. We also prove that the conjecture is true for prime $q \geq 3$, and indicate that more ingredients will be needed to settle Erd\"os’s conjecture for prime $q$.
\end{abstract}

\maketitle

\section{\bf Introduction}
\bigskip

In a written correspondence with Livingston, Erd\"os \cite{livingston} conjectured the following:
\begin{conjecture}(Erd\"os)\label{erdos-conjec}
Let $q$ be a positive integer and $f$ be an arithmetical function, periodic with period $q$. If $f(n) \in \{-1,1\}$ when $q \nmid n$ and $f(n) = 0$ otherwise, then 
\begin{equation*}\label{series-main}
\sum_{n=1}^{\infty} \frac{f(n)}{n} \neq 0,
\end{equation*}
whenever the series is convergent.
\end{conjecture}

In 1965, Livingston \cite{livingston} attempted to resolve the above conjecture. He predicted that to settle Conjecture \ref{erdos-conjec}, one would first have to prove: 
\begin{conjecture}(Livingston)\label{livingston-conjec}
Let $q \geq 3$ be a positive integer. The numbers 
\begin{equation*}
\bigg\{ \log\bigg(2 \sin \frac{a\pi}{q}\bigg) : 1 \leq a < \frac{q}{2} \bigg\} \text{  and } \pi
\end{equation*}
when $q$ is odd, and
\begin{equation*}
\bigg\{ \log\bigg(2 \sin \frac{a\pi}{q}\bigg) : 1 \leq a < \frac{q}{2} \bigg\} , \hspace{1mm} \pi, \text{ and } \log 2
\end{equation*}
when $q$ is even, are linearly independent over the field of algebraic numbers.
\end{conjecture}

The above statement does not depend on the branch of logarithm considered, as the values would only differ by an integer multiple of $2 \pi i$. In this paper, we disprove Livingston's conjecture in the case when $q$ is not prime and show that the conjecture is true when $q$ is prime. More precisely, we prove the following theorems:

\begin{theorem}\label{theorem-1}
Conjecture \ref{livingston-conjec} does not hold for $q \geq 4$ and $q$ not prime. In fact, for a composite positive integer $q \geq 6$, the numbers
\begin{equation*}
\bigg\{ \log \bigg( 2 \sin \frac{a \pi}{q} \bigg) : 1 \leq a < \frac{q}{2} \bigg\} 
\end{equation*}
are ${\mathbb{Q}}$-linearly dependent.
\end{theorem}

\begin{theorem}\label{theorem-2}
Let $p$ be an odd prime. The numbers 
\begin{equation*}
\bigg\{ \log \bigg( 2\sin \frac{a \pi}{p} \bigg) : 1 \leq a \leq \frac{p-1}{2} \bigg\} \hspace{2mm} \text{and } \pi
\end{equation*}
are $\bar{\mathbb{Q}}$-linearly independent. Thus, Conjecture \ref{livingston-conjec} is true when the modulus $p$ is an odd prime.
\end{theorem}

In both the above theorems, log denotes the principal branch. As a corollary of Theorem \ref{theorem-2}, we have
\begin{corollary}\label{coro}
Let $p$ be an odd prime and $f$ be an arithmetical function, periodic with period $p$ such that $f(n) \in \{-1,1\}$ when $p \nmid n$ and $f(n) = 0$ otherwise. Assume that $\sum_{a=1}^p f(a) = 0$. Then only one of the following is true, either
\begin{equation*}
\sum_{n=1}^{\infty} \frac{f(n)}{n} \neq 0,
\end{equation*}
or
\begin{equation*}
\sum_{a=1}^{p-1} f(a) \cot \bigg( \frac{a \pi}{p} \bigg) = \sum_{a=1}^{p-1} f(a) \cos \bigg( \frac{2 \pi ab}{p} \bigg) = 0,
\end{equation*}
for $1 \leq b \leq (p-1)/2$.
\end{corollary}

\section{\bf Preliminaries}
\bigskip

This section introduces some results that are fundamental to the proofs. 

\subsection{Baker's theorem on linear forms in logarithm of algebraic numbers}
We will use an important theorem of Baker (see \cite[Theorem 2.1, pg. 10]{baker}) concerning linear forms in logarithms of algebraic numbers, namely,
\begin{theorem} \label{baker-thm}
If $\alpha_1, \alpha_2, \cdots, \alpha_n$ are non-zero algebraic numbers such that $\log \alpha_1$, $\log \alpha_2$, $\cdots$, $\log \alpha_n$ are linearly independent over the rationals, then $1$, $\log \alpha_1$, $\log \alpha_2$, $\cdots$, $\log \alpha_n$ are linearly independent over the field of all algebraic numbers.
\end{theorem}

\subsection{Matrices of the Dedekind type}
Let $\mathfrak{M}$ be an $n \times n$ matrix with complex entries. Let $m_{i,j}$ denote the $(i,j)$-th entry of $\mathfrak{M}$. Then, $\mathfrak{M}$ is said to be of Dedekind type if there exists a finite abelian group, $G = \{ x_1, x_2, \cdots, x_n \}$ and a complex valued function $f$ on $G$ such that 
\begin{equation*}
m_{i,j} = f(x_i^{-1}x_j),
\end{equation*}
for all $1 \leq i,j \leq n$. We will use the following widely known theorem regarding matrices of the Dedekind type:
\begin{theorem}\label{Dedekind-det}
Let $\mathfrak{M}$ be an $n \times n$ matrix of the Dedekind type. For a character $\chi$ on $G$ ( a homomorphism of $G$ into $\mathbb{C}^*$), define
\begin{equation*}
S_{\chi} := \sum_{s \in G} f(s) \chi(s). 
\end{equation*}
Then the determinant of $\mathfrak{M}$ is equal to
\begin{equation*}
\prod_{\chi} S_{\chi},
\end{equation*}
where the product runs over all characters of $G$. Thus, $\mathfrak{M}$ is invertible if and only if
\begin{equation*}
S_{\chi} \neq 0,
\end{equation*}
for all characters $\chi$ of $G$.
\end{theorem}

For a proof of the above theorem and an exposition on properties of matrices of the Dedekind type, we refer the reader to \cite{ram-kaneenika}. The determinant of a matrix of the Dedekind type is often referred to as Dedekind determinant.

\subsection{Linear forms in logarithm of algebraic numbers with Dirichlet coefficients}
A Dirichlet character $\chi$ modulo $q$ is a group homomorphism,
\begin{equation*}
    \chi : {\big( \mathbb{Z}/q\mathbb{Z} \big)}^* \rightarrow \mathbb{C}^*,
\end{equation*}
which can be extended to a periodic function on all of integers by setting
\begin{equation*}
    \chi(n) = 
    \begin{cases}
    \chi(n \bmod q) & \text{if } (n,q) = 1, \\
    0 & \text{otherwise.}
    \end{cases}
\end{equation*}

The trivial Dirichlet character, $\chi_0$ is given by
\begin{equation*}
    \chi_0(n) = 
    \begin{cases}
    1 & \text{if } (n,q) = 1, \\
    0 & \text{otherwise.}
    \end{cases}
\end{equation*}

The Dirichlet $L$-function associated to a Dirichlet character $\chi$ is defined as
\begin{equation*}
    L(s,\chi) = \sum_{n=1}^{\infty} \frac{\chi(n)}{n^s},
\end{equation*}
which converges absolutely for $\Re(s) > 1$. The series $L(s,\chi)$ can be analytically continued to the entire complex plane except when $\chi = \chi_0$, in which case the series has a simple pole at $s = 1$. Since $\chi$ is a periodic arithmetical function, the proof of analytic continuation of $L(s,\chi)$ follows from the analytic continuation of the series $L(s,f)$ for a periodic arithmetical function $f$, proved in the next section, and the fact that $\sum_{a=1}^q \chi(a) = 0$ for a non-trivial Dirichlet character $\chi$ modulo $q$. We will make use of the following well-known lemma towards proving Theorem \ref{theorem-2}.

\begin{lemma}\label{l-1-chi}
Let $\chi$ be a non-trivial even Dirichlet character modulo an odd prime $p$, i.e, $\chi(-1) = 1$. Then,
\begin{equation*}
    \sum_{a=1}^{p-1} \overline{\chi}(a) \log  \big| 1 - \zeta_p^a \big| =  - \frac{p}{\tau({\chi})} L(1,\chi),
\end{equation*}
where 
\begin{equation}\label{gauss-sum}
    \tau(\chi) = \sum_{a=1}^p \chi(a) \zeta_p^a,
\end{equation}
is the Gauss sum associated to $\chi$ and $\zeta_p = e^{2 \pi i/p}$.
\end{lemma}

In the interest of completion, we include a proof of the above lemma.

\begin{proof}
Let $\chi$ be a non-trivial even Dirichlet character modulo an odd prime $p$. Let $\widehat{\chi}$ denote the discrete Fourier transform of $\chi$, given by
\begin{equation*}
    \widehat{\chi}(k) := \frac{1}{p} \sum_{a=1}^p \chi(a) \zeta_p^{-ak}.
\end{equation*}
This can be inverted using the identity
\begin{equation}\label{fourier-inversion}
\chi(n) = \sum_{k=1}^p \widehat{\chi}(k) \zeta_p^{kn}.
\end{equation}

Substituting expression \eqref{fourier-inversion} in the definition of the Dirichlet $L$-function associated to $\chi$ and noting that $\widehat{\chi}(p) = \sum_{a=1}^p \chi(a) = 0$ for a non-trivial Dirichlet character $\chi$, we get
\begin{align}
L(s,\chi) & = \sum_{n=1}^{\infty} \frac{1}{n^s} \sum_{k=1}^{p-1} \widehat{\chi}(k) \zeta_p^{kn}. \nonumber \\
 & = \sum_{k=1}^{p-1} \widehat{\chi}(k) \sum_{n=1}^{\infty} \frac{\zeta_p^{kn}}{n^s}. \label{inner-sum}
\end{align}

The inner sum converges for $s=1$. To see this, recall the partial summation formula,
\begin{theorem*}\label{abel-sum}
Let $\{a_n\}_{n \in \mathbb{N}}$ be a sequence of complex numbers and $f$ be a $C^{1}$ function on $\mathbb{R}_{>0}$. For $x > 0$, if $A(x):= \sum_{n \leq x} a_n$, then
\begin{equation*}
\sum_{1 \leq n \leq x} a_n f(n) = A(x)f(x) - \int_1^x A(t) f'(t) dt.
\end{equation*}
\end{theorem*}
For $1 \leq k \leq p-1$, let $a_n = \zeta_p^{kn}$ and $f(x) = 1/x$. Thus, $A(x) = \sum_{n \leq x} \zeta_p^{kn}$ and the partial summation formula gives us that
\begin{equation}\label{partial-sum}
\sum_{1 \leq n \leq x} \frac{\zeta_p^{kn}}{n} = \frac{A(x)}{x} + \int_1^x \frac{A(t)}{t^2} dt.
\end{equation}
Now, note that for $1 \leq k \leq p-1$,
\begin{equation*}
\sum_{n=1}^p \zeta_p^{kn} = 0.
\end{equation*}
Hence, the partial sums, $A(x)$ are bounded above by $p$ for all $x > 0$. Therefore, the integral in \eqref{partial-sum} is absolutely convergent as $x$ tends to infinity. Thus, taking limit as $x$ goes to infinity in \eqref{partial-sum}, we get the convergence of the inner sum in \eqref{inner-sum} and can conclude that
\begin{equation}\label{linear-form-log}
L(1,\chi) = - \sum_{k=1}^{p-1} \widehat{\chi}(k) \log(1 - \zeta_p^k),
\end{equation}
where $\log$ is the principal branch. Since $\chi$ is an even character, equation \eqref{linear-form-log} can be rewritten as
\begin{equation*}
\begin{split}
L(1,\chi) & = - \sum_{k=1}^{p-1} \widehat{\chi}(k) \log(1 - \zeta_p^k) \\
& = - \sum_{k=1}^{\floor*{(p-1)/2}} \widehat{\chi}(k) \big[ \log ( 1 - \zeta_p^k) + \log ( 1 - \zeta_p^{-k}) \big] \\
& = - \sum_{k=1}^{\floor*{(p-1)/2}} \widehat{\chi}(k) \log { \big| 1 - \zeta_p^k \big|}^2 \\
& = - \sum_{k=1}^{p-1} \widehat{\chi}(k) \log \big| 1 - \zeta_p^k \big|,
\end{split}
\end{equation*}
where $\widehat{\chi}$ denotes the Fourier transform of $\chi$. Now, note that the Fourier transform of $\chi$ can be evaluated in terms of the Gauss sum $\tau(\chi)$ as follows: for every $(k,p)=1$,
\begin{equation*}
\begin{split}
\widehat{\chi}(k) & = \frac{1}{p} \sum_{a=1}^{p-1} \chi(a) \zeta_p^{-ak} \\
& = \frac{1}{p} \sum_{t=1}^{p-1} \chi(-tk^{-1}) \zeta_p^t \\
& = \frac{\overline{\chi(-k)}}{p} \sum_{t=1}^{p-1} \chi(t) \zeta_p^t \\
& = \frac{\overline{\chi(-k)}}{p} \tau(\chi). 
\end{split}
\end{equation*}

Thus, the $L(s,\chi)$ for a non-trivial Dirichlet character $\chi$ has the value
\begin{equation*}\label{L(1,chi)}
L(1,\chi) = - \frac{\tau(\chi)}{p} \sum_{k=1}^{p} \overline{\chi}(k) \log | 1 - \zeta_p^k |
\end{equation*}
at $s = 1$. Another elementary but important fact about the Gauss sum is that when $\chi$ is a non-trivial Dirichlet character modulo $p$, 
\begin{equation}\label{gauss-sum-not-zero}
\tau(\chi) \neq 0.
\end{equation}
For a proof of the above fact, we refer the reader to \cite[Theorem 5.3.3, pg. 76]{murty}. This proves Lemma \ref{l-1-chi}.
\end{proof}

\section{\bf The approach of Livingston}
\bigskip

We first review general theory of $L$-series attached to a periodic arithmetical function following \cite{ram-saradha}. Let $q$ be a positive integer and $f$ be an arithmetical function that is periodic with period $q$. We define
\begin{equation*}
L(s,f) = \sum_{n=1}^{\infty} \frac{f(n)}{n^s}.
\end{equation*}

Let us observe that $L(s,f)$ converges absolutely for $\Re(s) > 1$. Since $f$ is periodic,
\begin{equation*}\label{hurw-form}
\begin{split}
L(s,f) & = \sum_{a=1}^q f(a)\sum_{k=0}^{\infty} \frac{1}{{(a + kq)}^s}\\
& = \frac{1}{q^s} \sum_{a=1}^q f(a) \zeta(s, a/q),
\end{split}
\end{equation*}
where $\zeta(s,x)$ is the Hurwitz zeta function. For $\Re(s) > 1$ and $0 < x \leq 1$, recall that the Hurwitz zeta function is defined as
\begin{equation*}
\zeta(s,x) = \sum_{n=0}^{\infty} \frac{1}{{(n+x)}^s}.
\end{equation*}

In 1882, Hurwitz \cite{hurwitz} proved that $\zeta(s,x)$ has an analytic continuation to the entire complex plane except for a simple pole at $s=1$ with residue $1$. In particular,
\begin{equation}\label{hurw}
\zeta(s,x) = \frac{1}{s-1} - \psi(x) + O(s-1),
\end{equation}
where $\psi$ is the digamma function, which is defined as the logarithmic derivative of the gamma function. This can be used to conclude that $L(s,f)$ can be extended analytically to the entire complex plane except for a simple pole at $ s=1 $ with residue $ \frac {1} {q} \sum_{a=1}^q f(a)$. Thus, $\sum_{n=1}^\infty \frac {f(n)} {n}$ exists if and only if $\sum_{a=1}^q f(a) = 0$, which we will assume henceforth. 

Let us also note that \eqref{hurw} helps us to express $L(1,f)$ as a linear combination of values of the digamma function. Therefore, 
\begin{equation}\label{digamma}
L(1,f) = - \frac{1}{q} \sum_{a=1}^{q} f(a) \psi \bigg( \frac{a}{q} \bigg).
\end{equation}

Let $f$ be an Erd\"os function, i.e, $f(n) = \pm 1$ when $q \nmid n$ and $f(n) = 0$ whenever $q|n$. The condition for the existence of $L(1,f)$ implies that
\begin{equation}\label{condition}
\sum_{a=1}^q f(a) = \sum_{a=1}^{q-1} f(a) = 0.
\end{equation} 

As seen earlier, $L(1,f)$ can be written as a linear combination of the values of the digamma function. Gauss (\cite[pg. 35-36]{davis}) proved the following formula for $1 \leq a < q$:
\begin{multline}\label{gauss-1}
\psi\bigg(\frac{a}{q} \bigg) = -\gamma - \log q - \frac{\pi}{2} \cot \bigg( \frac{a \pi}{q} \bigg) \\
+ \sum_{b=1}^{r} \bigg\{ \cos \bigg( \frac{2 \pi ab}{q} \bigg) \log \bigg( 4 \sin^2  \frac{\pi b}{q} \bigg) \bigg\} + {(-1)}^a \log 2 \hspace{1mm} \frac{1+{(-1)}^q}{2},
\end{multline} 
where $r := \floor*{(q-1)/2}$.

Substituting \eqref{gauss-1} in \eqref{digamma}, we have
\begin{multline*}
L(1,f) = \frac{-1}{q} \bigg[ \sum_{a=1}^{q-1} f(a) \bigg\{ \gamma + \log q + \frac{\pi}{2} \cot \bigg( \frac{a \pi}{q} \bigg) - \\
\sum_{b=1}^{r} \bigg\{ \cos \bigg( \frac{2 \pi ab}{q} \bigg) \log \bigg( 4 \sin^2 \frac{\pi b}{q} \bigg) \bigg\} + {(-1)}^a \log 2 \hspace{1mm} \frac{1+{(-1)}^q}{2} \bigg\} \bigg].
\end{multline*}

On simplifying the above expression using \eqref{condition}, we get
\begin{multline} \label{livingston-cdn}
L(1,f) =  \frac{-\pi}{2q} \sum_{a=1}^{q-1} f(a) \cot \bigg( \frac{a \pi}{q} \bigg)  \\
+ \frac{2}{q} \sum_{b=1}^{r} \bigg\{ \bigg[ \sum_{a=1}^{q-1} f(a) \cos \bigg( \frac{2 \pi ab}{q} \bigg) \bigg] \log \bigg( 2 \sin  \frac{\pi b}{q} \bigg)\bigg\} -  T_q,
\end{multline}
where
\begin{equation*}
T_q =
\begin{cases}
\frac{\log 2}{q} \bigg( \sum\limits_{k=1}^{q-1} {(-1)}^k f(k) \bigg) & \text{if $q$ is even}\\
0 & \text{otherwise.}
\end{cases}
\end{equation*}

Let us note that the numbers
\begin{equation*}
\cot \bigg( \frac{a \pi}{q} \bigg) \hspace{2mm} \text{and} \hspace{2mm} \cos \bigg( \frac{2 \pi ab}{q} \bigg)
\end{equation*}
are algebraic for $1 \leq a < q$ and $1 \leq b < q$. Since $f(a) \in \bar{\mathbb{Q}}$ and $f(q) = 0$, we are led to deduce that $L(1,f)$ is an algebraic linear combination of
\begin{equation*}
\pi, \log\bigg(2 \sin \frac{\pi}{q}\bigg), \log\bigg(2 \sin \frac{2\pi}{q}\bigg), \cdots, \log\bigg(2 \sin \frac{(q-1)\pi}{2q}\bigg)
\end{equation*}
together with $\log(2)$ when $q$ is even. This led Livingston to predict that if Conjecture \ref{erdos-conjec} were to be true, the above numbers should be linearly independent over $\bar{\mathbb{Q}}$. At this point, we make the following key observation - to conclude Conjecture \ref{erdos-conjec} as an implication of Conjecture \ref{livingston-conjec}, one is still required to prove that the resulting relation is non-trivial. That is, if $f$ is an Erd\"os function, not identically zero, then at least one of
\begin{equation}\label{cot-eqn}
\sum_{a=1}^{q-1} f(a) \cot \bigg( \frac{a \pi}{q} \bigg),
\end{equation}
or
\begin{equation}\label{cos-eqn}
\sum_{a=1}^{q-1} f(a) \cos \bigg( \frac{2 \pi ab}{q} \bigg), \hspace{2mm}  1 \leq b \leq r
\end{equation}
or $T_q$ is not zero. This question is not addressed by Conjecture \ref{livingston-conjec} and hence, Livingston's conjecture alone is not sufficient to settle the conjecture of Erd\"os.

\begin{remark*}
If $f$ is allowed to take values in $\bar{\mathbb{Q}}$ and $q$ is odd, then there exist a plethora of examples of functions $f$ that are not identically zero but for which \eqref{cot-eqn} and \eqref{cos-eqn} are both zero for all $1 \leq b \leq r$. These are given by the following theorem from \cite{bbw}:
\begin{theorem*}\label{bbw-odd}
Let $q \geq 3$ be a natural number. Then all odd, algebraically-valued functions $f$, periodic mod $q$, for which $L(1,f) = 0$ are given by the totality of  linear combinations with algebraic coefficients of the following $\floor*{\frac{1}{2}(q-3)}$ functions:
\begin{equation}
f_l(n) = ({-1})^{n-1} {\bigg( \frac{\sin n \pi/q}{\sin \pi/q} \bigg)}^l, \hspace{5mm} \text{for} \hspace{1mm} l=3,5,\cdots,(q-2) 
\end{equation}
when $q$ is odd and
\begin{equation*}
f_l(n) = ({-1})^{n-1} \bigg( \frac{\cos n \pi/q}{\cos \pi/q} \bigg) {\bigg( \frac{\sin n \pi/q}{\sin \pi/q} \bigg)}^l \hspace{3mm} \text{for} \hspace{1mm} l=3,5,\cdots,(q-1)
\end{equation*}
when $q$ is even. The functions are linearly independent and take values in $\mathbb{Q}(\zeta_q)$, i.e, the $q$-th cyclotomic field.
\end{theorem*}

Each $f_l$ in the above theorem is an odd function. Since $\cos(2 \pi ab/q)$ is an even function for $1 \leq a < q$, \eqref{cos-eqn} is zero for all $1 \leq b \leq r$. $T_q = 0$ as $q$ is odd. Thus, 
\begin{equation*}
L(1,f) = \frac{-\pi}{2q} \sum_{a=1}^{q-1} f(a) \cot \bigg( \frac{a \pi}{q} \bigg),
\end{equation*}
which is zero by the above theorem from \cite{bbw}.
\end{remark*}

\section{\bf Proof of the main theorems}
\bigskip

We make a useful observation before proceeding with the proofs. If $q$ is a positive integer and $1 \leq a < q/2$, then
\begin{equation}\label{identity-sin}
2 \sin \frac{a \pi}{q} = \frac{e^{i a \pi/q } - e^{-i a \pi/q}}{i} = i e^{-i a \pi/q} ( 1 - \zeta_q^a),
\end{equation}
where $\zeta_q = e^{2 \pi i/q}$. Since 
\begin{equation*}
\sin \frac{a \pi}{q} > 0,
\end{equation*}
for $1 \leq a < q/2$ and log denotes the principal branch,
\begin{multline}\label{transformation}
\log {\bigg( 2 \sin\frac{a \pi}{q} \bigg)} = \log {\bigg( \big| 1 - \zeta_q^a \big| \bigg)} + i0 = \log {\bigg( \big| 1 - \zeta_q^a \big| \bigg)} \\
= \log {\bigg( \big| 1 - \zeta_q^{-a} \big| \bigg)} = \log {\bigg( 2 \sin\frac{(q-a) \pi}{q} \bigg)}.
\end{multline}

\subsection{Proof of Theorem \ref{theorem-1}}
Conjecture \ref{livingston-conjec} does not hold for $q = 4$ because the numbers in consideration, namely
\begin{equation*}
\log \bigg( 2 \sin \frac{\pi}{4} \bigg) = \log \sqrt{2} = \frac{1}{2} \log 2, \hspace{1mm} \log 2 \text{  and } \pi 
\end{equation*}
are $\mathbb{Q}$ - linearly dependent.

Henceforth, assume that $q \geq 6$. We prove the linear dependence of the numbers 
\begin{equation*}
\bigg\{ \log \bigg( 2 \sin \frac{a \pi}{q} \bigg) : 1 \leq a < \frac{q}{2} \bigg\} 
\end{equation*}
by giving an explicit $\mathbb{Q}$-relation among them. Before proceeding, we note that by \eqref{transformation}, it suffices to exhibit a relation among logarithms of cyclotomic numbers. Now, since $q$ is not prime, there is a divisor $d$ of $q$ such that $d \neq 1,q$. For such a divisor $d$, we have the following polynomial identity in $\mathbb{C}[X,Y]$:
\begin{equation*}
X^{q/d} - Y^{q/d} = \prod_{j=1}^{q/d} \big(X - \zeta_{q/d}^j Y \big),
\end{equation*}
where $\zeta_{q/d} =e^{2 \pi i d/q}$. Substituting $X=1$ and $Y = \zeta_q^a$ for $(a,q)=1$, we have
\begin{equation*}
1 - e^{2 \pi i a/d} = \prod_{j=1}^{q/d} \big( 1 - e^{2 \pi i (dj/q + a/q)} \big) = \prod_{j=1}^{q/d} \big( 1 - e^{2 \pi i (a + dj)/q} \big)
\end{equation*}

Thus, taking absolute values of both sides of the above equation gives us
\begin{equation*}
{\bigg( \big| 1 - \zeta_q^{aq/d} \big|\bigg)} = \prod_{j=1}^{q/d} {\bigg( \big| 1 - \zeta_q^{(a + dj)} \big| \bigg)}.
\end{equation*}

Taking logarithms of both sides, we obtain the following $\mathbb{Q}$-linear relation
\begin{equation*}
\log \bigg( \big| 1 - \zeta_q^{aq/d} \big| \bigg) - \sum_{j=1}^{q/d} \log \bigg( \big| 1 - \zeta_q^{(a + dj)} \big| \bigg) = 0,
\end{equation*}
for all $1 \leq a < q$ and $(a,q)=1$ and $d|q$, $d \neq 1,q$. Hence, using \eqref{transformation}, we have 
\begin{equation}\label{unchanged-reln}
\log {\bigg( 2 \sin \bigg( \frac{aq}{d} \hspace{1mm }\frac{\pi}{q} \bigg) \bigg)} - \sum_{j=1}^{q/d} \log \bigg( 2 \sin \frac{(a + dj) \pi}{q} \bigg) = 0.
\end{equation}

Since we want a linear relation among 
\begin{equation*}
\bigg\{ \log \bigg( 2 \sin \frac{a \pi}{q} \bigg) : 1 \leq a < \frac{q}{2} \bigg\}, 
\end{equation*} 
we will replace  $\log ( 2 \sin (b \pi/q) )$ by $\log ( 2 \sin ((q-b) \pi/q) )$ whenever $b \geq q/2$. This is valid by \eqref{transformation}. Now, we make the following observations. Suppose that there exists an integer $k$ such that $ 1 \leq k < q/2$ and 
\begin{equation*}
k \equiv a + dj \equiv a + dl \bmod q,
\end{equation*}
for some $1 \leq j,l \leq q/d$ and $j \neq l$. This implies that $q | d(j - l)$, which is impossible since $(j - l) < q/d$. Thus,
\begin{equation}\label{imp-cdn-1}
a + dj \not\equiv a + dl \bmod q,
\end{equation}
for $1 \leq j,l \leq q/d$ and $j \neq l$. Similarly, 
\begin{equation}\label{imp-cdn-2}
-(a + dj) \not\equiv - (a + dl) \bmod q,
\end{equation}
for $1 \leq j,l \leq q/d$ and $j \neq l$. Suppose there exists a $k$ such that $ 1 \leq k < q/2$ and
\begin{equation*}
k \equiv a + dj \equiv - (a + dl) \bmod q,
\end{equation*}
for $1 \leq j,l \leq q/d$ and $j \neq l$. Thus, $q | (2a + d (j + l))$. Since $d | q$, we have $d | (2a + d (j - l))$, i.e, $d | 2a$. But $(a,q) = 1$. Hence, $(a,d) = 1$, which implies that $d | 2$. We assumed that $d \neq 1,q$. Therefore, $d = 2$. As a result, we have
\begin{equation}\label{imp-cdn-3}
a + dj \not\equiv - (a + dl) \bmod q,
\end{equation}
for $1 \leq j,l \leq q/d$ and $j \neq l$ unless $d = 2$. 

Thus, for $(a,q) = 1$, $d |q$ and $2 < d < q$, \eqref{unchanged-reln} along with \eqref{imp-cdn-1}, \eqref{imp-cdn-2} and \eqref{imp-cdn-3} give us a non-trivial $\mathbb{Q}$-relation, namely,
\begin{equation*}
\mathfrak{R}_{a,d}:= \sum_{1 \leq k < q/2} \alpha_k \log \bigg( 2 \sin \frac{k \pi}{q} \bigg) = 0,
\end{equation*}
where $\alpha_k$ is determined as follows: 
\begin{equation*}
\alpha_k = -1 \text{ if}
\begin{cases}
\text{either }  (aq/d \bmod q) < q/2, \hspace{1mm} k \not\equiv aq/d \bmod q \text{ \& } k \equiv \pm(a + dj) \bmod q \\
\text{or }  (aq/d \bmod q) \geq q/2, \hspace{1mm} k \not\equiv -(aq/d) \bmod q \text{ \& } k \equiv \pm (a + dj) \bmod q ,
\end{cases}
\end{equation*}
for some $1 \leq j \leq q/d$,
\begin{equation*}
\alpha_k = 1 \text{ if}
\begin{cases}
\text{either }  (aq/d \bmod q) < q/2, \hspace{1mm} k \equiv aq/d \bmod q \text{ \& } k \not\equiv \pm (a + dj) \bmod q \\
\text{or } (aq/d \bmod q) \geq q/2, \hspace{1mm} k \equiv  -(aq/d) \bmod q \text{ \& } k \not\equiv \pm (a + dj) \bmod q ,
\end{cases}
\end{equation*}
for some $1 \leq j \leq q/d$ and
\begin{equation*}
\alpha_k = 0, \text{ otherwise.}
\end{equation*}

To see that the above relation is non-trivial for $q$ not prime and $q \geq 6$, note that at least one of the following scenarios happens- either $(aq/d \bmod q) < q/2$, in which case for $ k \equiv aq/d \bmod q$, $\alpha_k = \pm 1$, or $(aq/d \bmod q) \geq q/2$, in which case for $k \equiv  -(aq/d) \bmod q$, $\alpha_k = \pm 1$.

Hence, the numbers under consideration in Conjecture \ref{livingston-conjec} are $\mathbb{Q}$-linearly dependent. As a result, Livingston's conjecture is false when $q$ is a composite number greater than or equal to $4$.

\subsection{Proof of Theorem \ref{theorem-2}}
We use the theory of Dedekind determinants developed in \cite{ram-kaneenika} and our knowledge of Dirichlet $L$-functions to prove that Conjecture \ref{livingston-conjec} is true when the modulus $q$ is prime. Consequently, let $p$ be an odd prime. Our aim is to prove that the numbers 
\begin{equation*}
\bigg\{ \log \bigg( 2\sin \frac{a \pi}{p}\bigg) : 1 \leq a \leq \frac{p-1}{2} \bigg\} \hspace{2mm} \text{and } \pi
\end{equation*}
are $\bar{\mathbb{Q}}$-linearly independent. 

Suppose, to the contrary, that the above numbers have a $\bar{\mathbb{Q}}$-linear relation among them. Thus, there exist algebraic numbers $\beta_0, \beta_1, \cdots, \beta_r$, not all zero, such that
\begin{equation}\label{relation}
\beta_0 \pi + \sum_{a=1}^r \beta_a \log \bigg(  2 \sin \frac{a \pi}{p} \bigg) = 0,
\end{equation}
where $r = (p-1)/2$. If $\beta_0 \neq 0$, then \eqref{relation} does not hold by the following Lemma from \cite{ram-saradha}:
\begin{lemma}\label{ram-saradha lemma}
If $c_0, c_1, \cdots, c_n$ are algebraic numbers and $\alpha_1, \alpha_2, \cdots, \alpha_n$ are positive algebraic numbers with $c_0 \neq 0$, then
\begin{equation*}
c_0 \pi + \sum_{j=1}^n c_j \log \alpha_j \neq 0.
\end{equation*}
\end{lemma}

Thus, $\beta_0$ must be zero. Now, if the numbers 
\begin{equation*}
\bigg\{ \log \bigg( 2 \sin \frac{a \pi}{p} \bigg) : 1 \leq a \leq \frac{p-1}{2} \bigg\} \hspace{2mm}
\end{equation*}
are $\mathbb{Q}$-linearly independent, then by Theorem \ref{baker-thm}, the above numbers are also $\bar{\mathbb{Q}}$-linearly independent. This contradicts our assumption, and hence, the above numbers must satisfy a $\mathbb{Q}$-linear relation. Thus, there exist $b_1, b_2, \cdots,b_r$ such that
\begin{equation}\label{Q-relation}
\sum_{a=1}^r b_a \log \bigg( 2 \sin \frac{a \pi}{p} \bigg) = 0.
\end{equation}

On clearing denominators, we can assume that
\begin{equation*}
b_a \in \mathbb{Z}, \hspace{1mm} 1 \leq a \leq \frac{(p-1)}{2}.
\end{equation*}

Since log denotes the principal branch and $\sin a \pi/p \in \mathbb{R}_{>0}$, \eqref{Q-relation} gives us the multiplicative relation
\begin{equation*}
\prod_{a=1}^r {\bigg( 2 \sin \frac{a \pi}{p} \bigg)}^{b_a} = 1.
\end{equation*}

Using \eqref{identity-sin}, this relation can be interpreted as a relation among roots of unity and cyclotomic numbers, i.e,
\begin{equation*}
\prod_{a=1}^r {\big( i e^{-i a \pi/p} ( 1 - \zeta_p^a)\big)}^{b_a} = 1.
\end{equation*}
The above relation can be further simplified by raising both sides of the equation to the $4p$-th power. Since ${(i e^{-ia \pi/p})}^{4p} = 1$, we are now left with the simpler multiplicative relation,
\begin{equation}\label{cycl-reln}
\prod_{a=1}^r {\big( 1 - \zeta_p^a \big)}^{B_a} = 1,
\end{equation}
where $B_a := 4p b_a$ and each factor in the product belongs to the cyclotomic field $\mathbb{Q}(\zeta_p)$.

Let $G$ be the group ${(\mathbb{Z}/p\mathbb{Z})}^{*} \big/ \{\pm1\}$. Let $c \in G$ and $\sigma_c$ be the unique automorphism of $\mathbb{Q}(\zeta_p)$ such that
\begin{equation*}
\sigma_c(\zeta_p) = \zeta_p^c.
\end{equation*}

The action of $\sigma_{c^{-1}}$ on \eqref{cycl-reln} gives us
\begin{equation*}
\prod_{a=1}^r { \big(1 - \zeta_p^{ac^{-1}}\big)}^{B_a} = 1.
\end{equation*}

On taking logarithm of the above equation, we obtain the relation
\begin{equation}\label{exp-1}
\sum_{a=1}^r B_a \log \bigg({2\sin \frac{ac^{-1} \pi}{p}\bigg)} = 0,
\end{equation}
for all $1 \leq a \leq r$ and $1 \leq c \leq r$.

Define an $r \times r$ matrix $\mathfrak{M}$ whose $(a,c)^{\text{th}}$ entry is 
\begin{equation*}
\log {\bigg( 2 \sin \frac{ac^{-1} \pi}{p} \bigg)}.
\end{equation*}

Thus, \eqref{exp-1} can be rewritten as a matrix equation, i.e,
\begin{equation*}
\mathfrak{M}v = 0,
\end{equation*}
where $v$ the $r \times 1$ column vector with the $a^{\text{th}}$-entry being $B_a$. Since \eqref{Q-relation} was a non-trivial relation, $v\neq0$. This is possible only if $\text{det }\mathfrak{M} = 0$.

Let $\mathfrak{M}^T$ denote the transpose of $\mathfrak{M}$. Notice that $\mathfrak{M}^T$ is a matrix of the Dedekind type with $\mathfrak{f} : G \rightarrow \mathbb{C}$ given by
\begin{equation*}
\mathfrak{f}(a) = \log {\bigg( 2 \sin \frac{a \pi}{p} \bigg)},
\end{equation*}
where $G$ is as defined above. As mentioned in Theorem \ref{Dedekind-det}, $\mathfrak{M}^T$ is invertible if and only if
\begin{equation*}
S_{\chi} := \sum_{a=1}^r \mathfrak{f}(a) \chi(a) \neq 0,
\end{equation*}
for all characters $\chi$ of the group $G$. Observe that all characters of the group $G$ are precisely the even Dirichlet characters modulo $p$. Thus, for a non-trivial even Dirichlet character $\chi$, we can use \eqref{transformation} to express $S_{\chi}$ as:
\begin{equation*}
\begin{split}
S_{\chi} & = \sum_{a=1}^r \chi(a) \log {\bigg( 2 \sin \frac{a \pi}{p} \bigg)} \\
& = \sum_{a=1}^r \chi(a)  \log \bigg( \big| 1 - \zeta_p^a \big| \bigg) \\
& = \frac{1}{2} \sum_{a=1}^{p-1} \chi(a) \log \bigg( \big| 1 - \zeta_p^a \big| \bigg) \\
& = -\frac{p}{2\tau(\bar{\chi})} L(1,\bar{\chi}),
\end{split}
\end{equation*}
where the last equality follows from Lemma \ref{l-1-chi}. By a famous theorem of Dirichlet,
\begin{equation*}
L(1, \bar{\chi}) \neq 0,
\end{equation*}
for non-trivial Dirichlet character $\chi$. Therefore, $S_{\chi} \neq 0$ when $\chi$ is a non-trivial character on $G$.

Now, let $\chi_0$ be the trivial character on $G$, i.e, $\chi_0$ is the trivial Dirichlet character modulo $p$. Then the factor $S_{\chi_0}$ is
\begin{equation*}
\begin{split}
S_{\chi_0} & = \sum_{a=1}^r \mathfrak{f}(a) \\
& = \sum_{a=1}^r \log {\bigg( 2 \sin \frac{a \pi}{p} \bigg)} \\ 
& = \sum_{a=1}^r \log { \bigg( \big| 1 - \zeta_p^a \big| \bigg) } \\
& = \frac{1}{2} \log \bigg( \prod_{a=1}^{p-1} \big| 1 - \zeta_p^a \big| \bigg) \\
& = \frac{1}{2} \log p \neq 0,
\end{split}
\end{equation*}
where the last equality can be derived by noting that
\begin{equation*}
\frac{1 - X^p}{1-X} = \sum_{j=0}^{p-1} X^j = \prod_{a=1}^{p-1} ( 1 - \zeta_p^aX),
\end{equation*}
substituting $X=1$ and taking absolute values of both sides. Thus, $S_{\chi_0} \neq 0$. 

Hence, $\mathfrak{M}^T$, and in turn, $\mathfrak{M}$ is invertible. Therefore $v = 0$, which is a contradiction. This proves the theorem. 

\subsection{Proof of Corollary \ref{coro}}
Suppose that 
\begin{equation*}
\sum_{n=1}^{\infty} \frac{f(n)}{n} = L(1,f) = 0.
\end{equation*}

From Theorem \ref{theorem-2}, we see that Conjecture \ref{livingston-conjec} is true when the period of $f$ is an odd prime, i.e, that the numbers 
\begin{equation*}
    \bigg\{ \log \bigg( 2\sin \frac{a \pi}{p} \bigg) : 1 \leq a \leq \frac{p-1}{2} \bigg\} \hspace{2mm} \text{and } \pi
\end{equation*}
are $\bar{\mathbb{Q}}$ - linearly independent. Thus, the relation obtained from \eqref{livingston-cdn}, namely, 
\begin{equation*}
0 =  \frac{-\pi}{2p} \sum_{a=1}^{p-1} f(a) \cot \bigg( \frac{a \pi}{p} \bigg) 
+ \frac{2}{p} \sum_{b=1}^{r} \bigg\{ \bigg[ \sum_{a=1}^{p-1} f(a) \cos \bigg( \frac{2 \pi ab}{p} \bigg) \bigg] \log \bigg( 2 \sin  \frac{\pi b}{p} \bigg)\bigg\} 
\end{equation*}
is a trivial relation. Therefore, the co-efficients of $\pi$ and $\log ( 2 \sin (b \pi/p))$ must all be zero. This proves the corollary.

\section*{Acknowledgments}
I would like to extend my gratitude to Prof. M. Ram Murty for bringing this conjecture to my notice and discussing the various nuances involved in the proof. I am also thankful for his guidance and help in writing the note. I am very much obliged to the referee for insightful comments on earlier versions of this paper.

\end{document}